\documentclass[microtype]{gtpart}

\usepackage{pinlabel}  
\usepackage{graphicx}  
\usepackage[all]{xy}
\usepackage[symbol]{footmisc}
\usepackage[mathscr]{eucal}

\title{Moments of a Length Function on the Boundary of a Hyperbolic Manfiold}

\author{Nicholas G. Vlamis
\footnote{The author acknowledges support from U.S. National Science Foundation grants DMS 1107452, 1107263, 1107367 "RNMS: GEometric structures And Representation varieties" (the GEAR Network).}
}

\givenname{Nicholas}
\surname{Vlamis}
\address{Department of Mathematics \\ Boston College \\ 140 Commonwealth Ave. \\ Chestnut Hill, MA 02467}
\email{nicholas.vlamis@bc.edu}
\urladdr{http://nickvlamis.com}


\arxivreference{math.GT/1303.6494}
\arxivpassword{d9mdk}

\volumenumber{}
\issuenumber{}
\publicationyear{}
\papernumber{}
\startpage{}
\endpage{}
\doi{}
\MR{}
\Zbl{}
\received{}
\revised{}
\accepted{}
\published{}
\publishedonline{}
\proposed{}
\seconded{}
\corresponding{}
\editor{}
\version{}

\newcommand{\bc}{\mathbb C}

\newcommand{\bh}{\mathbb H}

\newcommand{\bz}{\mathbb Z}

\newcommand{\br}{\mathbb R}
\newcommand{\bS}{\mathbb S}

\newcommand{\wt}{\widetilde}

\newcommand{\al}{\alpha}
\newcommand{\be}{\beta}

\newcommand{\vp}{\varphi}

\newcommand{\Li}{\mathrm{Li}_2}
\newcommand{\li}{\mathscr{D}}

\DeclareMathOperator{\sech}{sech}
\DeclareMathOperator{\arccoth}{arccoth}
\DeclareMathOperator{\Isom}{Isom}
\DeclareMathOperator{\Vol}{Vol}

\newtheorem{Thm}{Theorem}[section]
\newtheorem{Prop}[Thm]{Proposition}
\newtheorem{Lem}[Thm]{Lemma}
\newtheorem{Cor}[Thm]{Corollary}

\theoremstyle{definition}
\newtheorem{Def}[Thm]{Definition}

\theoremstyle{remark}

\begin{document}


\begin{abstract}
In this paper we will study the statistics of the unit geodesic flow normal to the boundary of a hyperbolic manifold with non-empty totally geodesic boundary. Viewing the  time it takes this flow to hit the boundary as a random variable, we derive a formula for its moments in terms of the orthospectrum.  The first moment gives the average time for the normal flow acting on the boundary to again reach the boundary, which we connect to Bridgeman's identity (in the surface case), and the zeroth moment recovers Basmajian's identity. Furthermore, we are able to give explicit formulae for the first moment in the surface case as well as for manifolds of odd dimension.   In dimension two, the summation terms are dilogarithms.  In dimension three, we are able to find the moment generating function for this length function.
\end{abstract}

\maketitle


\section{Introduction}

Let $M$ be a compact hyperbolic manifold with non-empty totally geodesic boundary.  An {\em orthogeodesic} for $M$ is an oriented geodesic arc with endpoints normal to $\partial M$ (see \cite{basmajian}). We will denote the collection of orthogeodesics by $O_M = \{\al_i\}$.   Let $\ell_i$ denote the length of $\al_i$, then the collection $|O_M| = \{\ell_i\}$ (with multiplicities) is known as the {\em orthospectrum}.  As we will be summing over the orthospectrum, it is important to note that $O_M$ is a countable collection: this can be seen by doubling the manifold and observing that the orthogeodesics correspond to a subset of the closed geodesics in the double.

Given $x\in \partial M$, let $\al_x$ be the geodesic emanating from $x$ normal to $\partial M$.  Then, as the limit set is measure zero, for almost every $x\in \partial M$ we have that $\al_x$ terminates in $\partial M$; hence, the length of $\al_x$ is finite.  This allows us to define the measurable function  $L\co \partial M \to \br$ given by $L(x) = length(\alpha_x)$.  Let $dV$ denote the hyperbolic volume measure on $\partial M$, then $V(\partial M)$ is finite allowing us to define the  probability measure $dm = dV/V(\partial M)$ on $\partial M$, so that $(\partial M, dm)$ is a probability space.  This lets us view $L\co\partial M \to \br$ as a random variable.  Given a random variable $X$ on a probability space with measure $p$, the {\it $k^{th}$-moment of $X$} is defined to be $E[X^k] = \int X^k dp$, where $E[X]$ denotes the expected value. Let $A_k(M)$ be the $k^{th}$ moment of $L$.  In particular, $A_1(M)$ is the expected value of $L$.  In this paper we will show that the positive moments of $L$ are finite and encoded in the orthospectrum:

\begin{Thm}\label{finiteness}
Let $M= M^{n}$ be an $n$-dimensional compact hyperbolic manifold with nonempty totally geodesic boundary, then $A_k(M)$ is finite for all $k\in \bz^{\geq 0}$.
\end{Thm}

\begin{Thm}\label{ortho}
Let $M= M^{n}$ be an $n$-dimensional compact hyperbolic manifold with nonempty totally geodesic boundary, then for all $k\in \bz^{\geq 0}$
$$A_k(M) = \frac{1}{V(\partial M)} \sum_{\ell \in |O_M|} F_{n,k}(\ell), $$
where $$F_{n,k}(x) =\Omega_{n-2}\int_0^{\log\coth(x/2)}\left[ \log\left(\frac{\coth x +\cosh r }{\coth x-\cosh r}\right)\right]^k \sinh^{n-2}(r)\ dr$$
and  $\Omega_n$ is the volume of the standard $n$-sphere.  Furthermore, the identity for $A_0(M)$ is Basmajian's identity.
\end{Thm}

\noindent
Basmajian's identity gives the volume of the boundary in terms of the orthospectrum:
\begin{Thm}[Basmajian's Identity, \cite{basmajian}]
If $M$ is a compact hyperbolic $n$-manifold with totally geodesic boundary, then
$$V(\partial M) = \sum_{\ell_i\in|O_M|} V_{n-1}\left(\log\coth\frac {\ell_i}2\right),$$
where $V_n(r)$ is the volume of the hyperbolic $n$-ball of radius $r$.
\end{Thm}
\noindent
Note that by combining Theorem \ref{ortho} and Basmajian's identity we see that $A_k(M)$ depends solely on the orthospectrum.

As corollaries to Theorem \ref{ortho} we can write the function $F_{n,1}(x)$ in dimension 2  and all odd dimensions without integrals.  In the following corollary $\Li(x)$ is the standard dilogarithm (see \cite{lewin}).  We will also write $\ell(\partial S)$ for sum of the lengths of each boundary component of a surface $S$.

\begin{Cor}\label{surface}
Let $S$ be a compact hyperbolic surface with nonempty totally geodesic boundary.  Then
$$A_1(S) = \frac2{\ell(\partial S)} \sum_{\ell\in |O_S|}\left[\Li\left(-\tanh^2\frac {\ell}2\right)-\Li\left(\tanh^2\frac {\ell}2\right)+\frac {\pi^2}4\right]$$
\end{Cor}

\begin{Cor}\label{odd}
Let $M$ be an $n$-dimensional compact hyperbolic manifold with nonempty totally geodesic boundary where $n$ is odd.  Then
$$\resizebox{1\hsize}{!} {$\displaystyle{A_1(M) = \frac{2\Omega_{n-2}}{Vol(\partial M)}\sum_{\ell\in |O_M|}\sum_{j=0}^{\frac{n-3}2}\frac{(-1)^{\frac{n-3}2-j}\binom{\frac{n-3}2}{j}}{2j+1}\coth^{2j+1}(\ell)\left[\log(2\cosh\ell)-\ell_i\tanh^{2j+1}(\ell)+\sum_{k=1}^j\frac{1-\tanh^{2k}(\ell)}{2k}\right].}$}$$
\end{Cor}

The rest of the paper is dedicated to understanding the asymptotics of the $F_{n,k}$'s and finding the moment generating function in dimension 3.
The motivation of this paper comes from recent work of Bridgeman and Tan in \cite{bt}, where the authors study the moments of the hitting function associated to the unit tangent bundle of a manifold (i.e. the time it takes the geodesic flow of a vector to reach the boundary).  In the paper they are able to show the moments are finite and give an explicit formula for the expected value in the surface case as well as relate the orthospectrum identities of Basmajian and Bridgeman (see \cite{b}, \cite{bk}, and \S6 below) as different moments of the hitting function.  In \S6 we give a relationship between Bridgeman's identity and $A_1(S)$ in dimension 2.

\subsection*{Acknowledgements}
The author would like to thank his advisor, Martin Bridgeman, for his guidance as well as Ser Peow Tan for helpful conversations and hosting the author at the National University of Singapore.  In addition, the author would like to thank the reviewer for valuable comments.


\section{Kleinian Groups}

For $n\geq 2$, let $\Isom^+(\bh^n)$ be the space of orientation preserving isometries of hyperbolic $n$-space.  With the topology of uniform convergence on the space of isometries, we define a {\it Kleinian group} to be a discrete torsion-free subgroup of $\Isom^+(\bh^n)$. If $\Gamma<\Isom^+(\bh^n)$ is a Kleinian group, then $\bh^n/\Gamma$ is a hyperbolic manifold, i.e. a Riemannian manifold of constant curvature $-1$.

In the Poincar\'e model of hyperbolic space we can identify the boundary of $\bh^n$ with the $(n-1)$-sphere called the {\it sphere at infinity} and denoted $S^{n-1}_\infty$.  Pick $x\in \bh^n$ and define the {\it limit set of $\Gamma$} to be the set $\Lambda_\Gamma = \overline{\Gamma x} \cap S^{n-1}_\infty$.  Note that this definition is independent of the choice of $x$.  Define the {\it convex hull} $\mathrm{CH}(\Lambda_\Gamma)$ of the limit set $\Lambda_\Gamma$ to be the smallest convex subset  of $\bh^n$ containing all the geodesics in $\bh^n$ with endpoints in $\Lambda_\Gamma$. As $\Lambda_\Gamma$ is $\Gamma$-invariant, so is $\mathrm{CH}(\Lambda_\Gamma)$ and so we can take the quotient of $\mathrm{CH}(\Lambda_\Gamma)$ by $\Gamma$, which we call the {\it convex core} and denote $\mathrm{C}(\Gamma)$.  A Kleinian group is {\it convex cocompact} if its associated convex core is compact (see \cite{thurston}).

With these definitions at hand, we recall the following two theorems from Patterson-Sullivan theory (see \cite{nicholls}):

\begin{Thm}\label{thm:orbits}
Let $\Gamma<\mathrm{Isom}^+(\bh^n)$ be a convex cocompact Kleinian group and let $\delta = \delta(\Gamma)$ be the Hausdorff dimension of the limit set of $\Gamma$.  There exists $r_0$ such that  for $r\geq r_0$, $$N_x(r) = \left|\{\gamma\in \Gamma\colon d(\gamma(x),x)<r\}\right| \leq ae^{\delta r},$$
for some constant $a$ depending on $\Gamma$ and $x$.
\end{Thm}

\begin{Thm}
Let $\Gamma<\mathrm{Isom}^+(\bh^n)$ be a convex cocompact Kleinian group and let $\delta = \delta(\Gamma)$ be the Hausdorff dimension of the limit set of $\Gamma$. Then $\delta = n-1$ if and only if $\bh^n/\Gamma$ is finite volume. 
\end{Thm}


\section{Finite Moments}

Let $M = M^n$ be a compact $n$-dimensional hyperbolic manifold with totally geodesic boundary.  As above, let $L$ denote the time to the boundary of the unit normal flow on the boundary.  We let $dV$ be the induced hyperbolic volume measure on the boundary and define $dm = dV/V(\partial M)$, so that $(\partial M, m)$ is a probability space and $L:\partial M \to \br$ is a random variable on this space.  We let $A_k(\partial M) = E[L^k]=\int_{\partial M} L^k\, dm$ be the $k^{th}$ moment of $L$. In this section we will show that $A_k(M)$ is finite for all nonnegative integers $k$. We first state a basic fact from hyperbolic geometry.

\begin{Lem} \label{lem: projection}
Let $U$ be a hyperplane in $\bh^n$ and $B_R$ a hyperbolic $n$-ball of radius $R$ a distance $s$ from $U$.  The orthogonal projection of $B_R$ to $U$ has radius $r<\log\coth(s/2)$.
\end{Lem}

\begin{proof}
Let $p\in \partial B_R$ be the point closest to $U$, so that $d(p,U) = s$ and let $V$ be the hyperplane containing $p$ such that $d(U,V) = s$.  Then the orthogonal projection of $B_R$ is contained in the orthogonal projection of $V$.  The orthogonal projection of $V$ to $U$ has radius $\log\coth(s/2)$ (see \cite{basmajian}), implying that $r < \log\coth(s/2)$ as desired.
\end{proof}

We can now show that $A_k(M)$ is finite:

\medskip
\noindent
{\bf Theorem \ref{finiteness}} {\it
Let $M= M^{n}$ be an $n$-dimensional compact hyperbolic manifold with nonempty totally geodesic boundary, then $A_k(M)$ is finite for all $k\in \bz^{\geq 0}$.}
\begin{proof}
We want to work in hyperbolic space: identify the universal cover $\wt M$ of $M$ with a convex subset of $\bh^n$, so that $\pi_1M=\Gamma < \mathrm{Isom}^+(\bh^n)$ is a convex cocompact Kleinian group.  As $M$ has a finite number of disjoint boundary components and we are investigating the integral over the boundary, it is enough to prove finiteness for a single component.  Fix $S\subset \partial M$ a component and a lift $\wt S\subset \wt M$ of $S$ (note: $\wt S$ is a copy of $\bh^{n-1}$ sitting in $\bh^n$).  Let $U$ be a convex fundamental domain for the action of $\Gamma$ on $\wt M$. Pick $p\in U$ and let $B_R(p)$ be a ball centered at $p$ of radius $R$ such that $U\subset B_R(p)$.  Set $W= U\cap \wt S$ to be a fundamental domain for the action of $\mathrm{Stab}(\wt S)<\Gamma$ on $\wt S$. Define $n_t\co W \to \bh^n$ to be the unit geodesic flow normal to $\wt S$ into $\wt M$ for a time $t$ and set $X_t = \{x\in W\colon n_t(x) \in \wt M\}.$ Define $\pi\co \bh^n \to \wt S$ to be orthogonal projection.  

We will now bound $V(X_t)$ for $t\geq r_0$, where $r_0$ is taken from Theorem \ref{thm:orbits}.  If $x\in X_t$, then $n_t(x)\in \gamma U$ for some $\gamma\in \Gamma$. If $n_t(x)\in n_t(X_t)\cap \gamma U$, then $d(p,\gamma(p)) < t+2R$.  Let $\Gamma_t =\{\gamma\in \Gamma\colon n_t( X_t)\cap \gamma U \neq \emptyset\}$, then from the above theorem, we know that $|\Gamma_t| \leq N_p(t+2R) \leq ae^{\delta(t+2R)}$, where $\delta$ is the Hausdorff dimension of the limit set of $\Gamma$.  As $n_t(X_t)\subset \bigcup_{\gamma\in\Gamma_t} \gamma U$ and $\pi(n_t(X_t))=X_t$, we have
$$V(X_t) \leq \sum_{\gamma\in\Gamma_t}V(\pi(\gamma U)).$$
Now, fix $\gamma\in \Gamma_t$, then $\gamma\cdot U \subset B_R(\gamma\cdot p).$  Suppose that $B_R(\gamma\cdot p)$ is a distance $s$ from $W$ and let $r$ be the radius of its projection, we then have that $t < r+s+2R$ by the triangle inequality; in particular, $s>t-r-2R$.  Furthermore, as orthogonal projection is always distance decreasing in hyperbolic space, $r< R$, so that $s>t-3R$.  From the above lemma, we can conclude that $$r \leq \log\coth(s/2)\leq \log\coth\left(\frac{t-3R}2\right)\equiv f(t).$$
As the above bound for the radius does not depend on $\gamma$, we have 
$$V(X_t)\leq |\Gamma_t|V_{n-1}\left(f(t)\right)\leq N_p(t+2R) V_{n-1}\left(f(t))\right)\leq ae^{\delta(t+2R)}\,V_{n-1}\left(f(t)\right),$$
where $V_{n}(r)$ is the volume of a $n$-dimensional hyperbolic ball of radius $r$.  We observe two asymptotics:  1) $\lim_{x\to \infty} e^{x}\log\coth(x/2) = 2$ and 2) $\lim_{r\to0} V_n(r)/r^n  = C_n$ for some constant $C_n>0$.  From these facts and the above inequality, we see that $$\lim_{t\to\infty} e^{(n-1-\delta)t}\cdot V(X_t)\leq A,$$ for some constant $A$.  From the theorem stated above, we know that $n-1-\delta>0$.  

We now move to the moments.  We have setup the following situation:
\begin{align*}
\int_S L^k dV &= \sum_{t=0}^\infty \int_{L^{-1}(t,t+1)} L^k dV \leq \sum_{t=0}^{\infty} (t+1)^k\int_{L^{-1}(t,t+1)}dV
 \leq \sum_{t=0}^\infty (t+1)^k V(X_t).
\end{align*}
But, we saw that the asymptotics of $V(X_t)$ are less than a multiple of $e^{-bt}$ with $b>0$, which implies the above sum converges since $\sum(t+1)^ke^{-bt}$ converges.
\end{proof}


\section{The Moments as a Sum over the Orthospectrum}
\subsection{Basmajian's Ball Decomposition of the Boundary}

In this section we introduce a decomposition of $\partial M$ into a disjoint union of $n-1$ balls (affectionately known as ``leopard spots").  We will recall Danny Calegari's method of accomplishing this in \cite{calegari}.

\begin{Def}
Let $X$ and $Y$ be totally geodesic copies of $\bh^{n-1}$ sitting inside of $\bh^n$ with disjoint closure in $\bh^n\cup S^{n-1}_\infty$.  A {\em chimney} is the closure of the union of the geodesic arcs from $X$ to $Y$ that are perpendicular to $X$.  
\end{Def}

The distance between the hyperplanes $X$ and $Y$ defining the chimney is realized by a unique geodesic  perpendicular to both planes, called the {\em core}, the length of which is the {\em height} of the chimney. The chimney cuts out a disk in $X$, which is called the {\em base}.  Let $\al$ be the geodesic containing the core and $\be$ a geodesic containing a diameter of the base.  Then $\al$ and $\be$ span a copy of $\bh^2$ in $\bh^n$.  Furthermore, the intersection of this plane with the chimney cuts out half an ideal quadrilateral with orthospectrum $\{2\ell, 2r\}$, where $\ell$ is the length of the core and $r$ the radius of the base.  We then have $\sinh(r)\sinh(\ell) = 1$, which yields $r = \log\coth\frac{\ell}2$ (see \cite{beardon}). 

\begin{Thm}[Chimney Decomposition, \cite{calegari}]
Let $M$ be a compact hyperbolic $n$-manifold with totally geodesic boundary $S$.  Let $M_S$ be the covering space of $M$ associated to $S$.  Then $M_S$ has a canonical decomposition into a piece of zero measure, together with two chimneys of height $\ell_i$ for each $\ell_i$ in the (unordered) orthospectrum.
\end{Thm}

If we take the bases of the chimneys in the decomposition of the above theorem, then we get a decomposition of $\partial M$ into $(n-1)$-balls. With this decomposition, we can give the quick proof of Basmajian's identity in \cite{calegari}:
\begin{proof}[Proof of Basmajian's identity]
Recall that we are working with the ordered orthospectrum. $S$ in $M_S$ is decomposed into a set of measure zero together with the union of the bases of the chimneys.  Thus
$$V(S) = \sum_i V_{n-1}\left(\log\coth\frac{\ell_i}2\right),$$
where $V_n(r)$ is the volume of a hyperbolic $n$-ball of radius $r$.
\end{proof}


\subsection{Deriving the Length Function}

\begin{figure}[ht!]
\labellist
\small\hair 2pt
\pinlabel $\ell_i$ at 210 -12
\pinlabel $r$ at -10 160
\pinlabel {$\wt x$} at -5 305
\pinlabel $L_i(r)$ at 185 325
\endlabellist
\centering
\includegraphics[scale=0.5, clip=false]{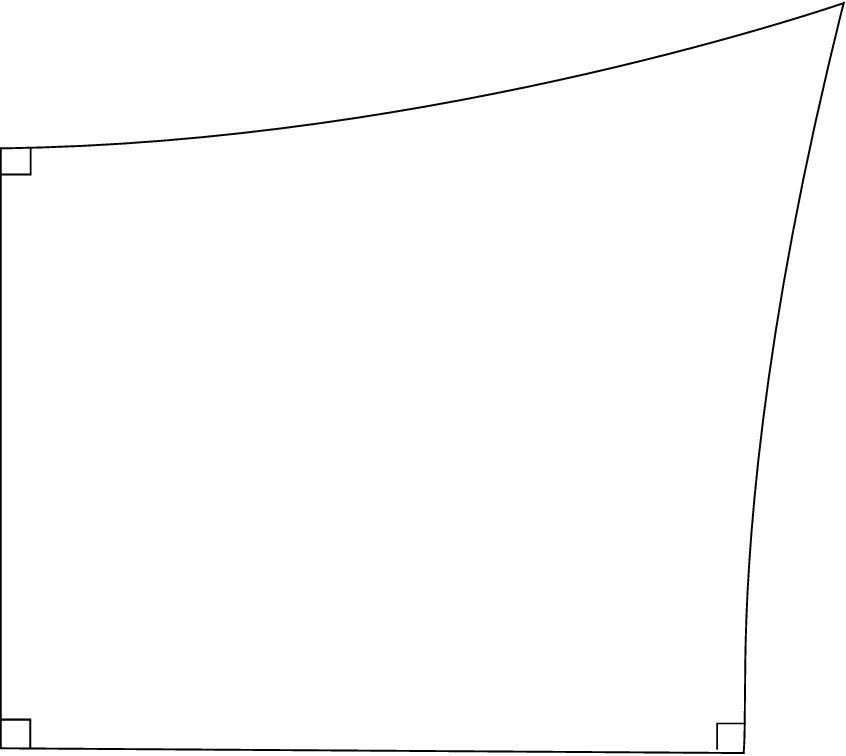}
\caption{A Lambert quadrilateral showing the setup for $L(r)$.}
\label{fig:figure}
\end{figure}
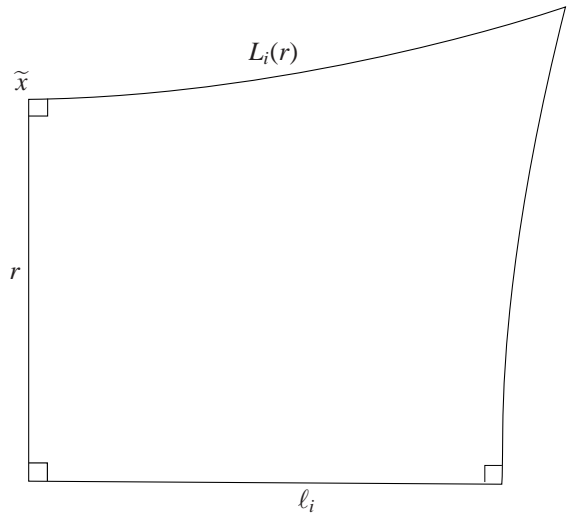

Let $U_i$ be the  interior of the ball associated to $\ell_i \in |O_M|$.  By above, the union of the $U_i$'s is a full measure set in $S = \partial M$.  The measurable function $L: S \to \br$ assigning to each $x\in S$ the length of the arc emanating perpendicularly from $S$ at $x$ can be written as $L = \sum_{\ell_i\in |O_M|} L_i$, where $L_i = L|U_i \co U_i \to\br$ since the $U_i$'s are pairwise disjoint. As a chimney has rotational symmetry about its core, we see that $L(x)$ for $x\in S$ depends only on the distance between $x$ and the core, ie $L_i$ is a function of the radius; hence, deriving a formula for $L_i$ is a problem in the hyperbolic plane.  Associated to each $U_i$ are two components of the boundary, $R_i$ and $T_i$,  and two lifts of these components to hyperplanes in $\bh^n$, $\wt R_i$ and $\wt T_i$.  If $x\in R_i$, then we are interested in the chimney with its base in $\wt R_i$ and the lift of $x$ sitting in $\wt R_i$, call it $\wt x$.  There is a unique copy of $\bh^2\subset \bh^n$ containing the core of the chimney, $\wt x$, and the geodesic connecting the two.  The geodesic contained in this plane going through $\wt x$ and intersecting $\wt R_i$ perpendicularly intersects $\wt T_i$;  furthermore, the length of this arc is $L_i(x)$. The diagram showing this situation in $\bh^2$ is shown in Figure \ref{fig:figure}.

\noindent
We see that $L_i(x)$ is the length of a side in a Lambert quadrilateral (a hyperbolic quadrilateral with three right angles). Let $r$ be the distance from $x$ to the core, then as we noted $L_i$ is solely a function of the radius, we will write $L_i(x) = L_i(r)$.  From hyperbolic trigonometry we find $\coth L_i(r) = \sech(r)\coth(\ell_i)$ (see \cite{beardon}) or 
\begin{equation} \label{eq: length}
 L_i(r) = \arccoth(\sech(r)\coth(\ell_i)) = \frac12 \log \left(\frac{\coth \ell_i + \cosh r}{\coth \ell_i - \cosh r}\right),
\end{equation}
where the second equality holds as $\sech(r)\coth(\ell_i) > 1$ on the domain of interest $r\in [0, \log\coth(\ell_i/2))$.


\subsection{Proof of Theorem 1.2}

For completeness, we restate the result:

\medskip
\noindent
{\bf Theorem \ref{ortho}.} 
{\it Let $M= M^{n}$ be an $n$-dimensional compact hyperbolic manifold with nonempty totally geodesic boundary, then for all $k\in \bz^{\geq 0}$
$$A_k(M) = \frac{1}{V(\partial M)} \sum_{\ell \in |O_M|} F_{n,k}(\ell), $$
where $$F_{n,k}(x) =\Omega_{n-2}\int_0^{\log\coth(x/2)}\left[ \log\left(\frac{\coth x +\cosh r }{\coth x-\cosh r}\right)\right]^k \sinh^{n-2}(r)\ dr$$
and  $\Omega_n$ is the volume of the standard $n$-sphere.  Furthermore, the identity for $A_0(M)$ is Basmajian's identity.}

\begin{proof}
From the additivity property of measures we have $\int L^k\,\, dm = \sum \int_{U_i} L_i^k \, dm$. As $dm = dV/V(\partial M)$ and $dV$ is the $(n-1)$-dimensional hyperbolic volume form we can write it in spherical coordinates as $$dm = \frac{1}{V(\partial M)}\sinh^{n-2}(r)\, dr\, d\Omega_{n-2},$$ where $d\Omega_{n-2}$ is the volume form on the standard unit sphere. Above we saw that $L_i$ is a function solely of the radius and $U_i$ is a $(n-1)$-dimensional hyperbolic ball with radius $\log(\coth(\ell_i/2))$, so that
\begin{align*}
\int_{U_i} L_i^k \, dm &= \frac1{V(\partial M)} \int_{S^{n-2}}\int_0^{\log(\coth(\ell_i/2))} L_i^k(r) \sinh^{n-2}(r)\, dr\, d\Omega_{n-2} \\
&= \frac{\Omega_{n-2}}{V(\partial M)} \int_0^{\log(\coth(\ell_i/2))} L_i^k(r) \sinh^{n-2}(r)\, dr,
\end{align*}
where we write $\Omega_{n-2}$ for the volume of the standard $(n-2)$-dimensional unit sphere. Define $F_{n,k}(x)$ as stated in the theorem, so that the equality holds for $A_k(M)$ by \eqref{eq: length}.

Observe that $F_{n,0}(x)$ gives the volume of a hyperbolic $(n-1)$-ball of radius $\log\coth(x/2)$.  As $A_0(M) = 1$, we see that this identity yields
$$1 = \frac{1}{V(\partial M)} \sum_{\ell\in |O_M|} V_{n-1}(\log\coth(\ell/2)),$$
which is Basmajian's identity.
\end{proof}


\section{Surface Case}
\subsection{Dilogarithms}

For $|z|<1$ in $\bc$ the dilogarithm is defined as $$\Li(z) = \sum_{n=1}^\infty \frac{z^n}{n^2}.$$  Using the Taylor series for $\log(1-z)$ about $z=0$, we can write
$$\Li(z) = \int_z^0 \frac{\log(1-z)}{z}dz.$$  One can then take a branch of $\log(z)$ in order to analytically continue $\Li(z)$ to the complex plane minus a branch cut. The standard definition of the dilogarithm assumes the branch cut for $\log(z)$ is along $(-\infty,0]$; however, for our purposes we will be interested in a different branch cut.  Define the function $\li(z)$ to be the resulting dilogarithm by using the branch cut along $(-i\infty,0]$ for $\log(z)$  such that $\log(-1) = i\pi$.  We note that $\Li(z) = \li(z)$ for $z\in (-\infty, 1)$. 

The dilogarithm $\Li(z)$ has the following well-known identity (see \cite{lewin}):
$$ \Li(z) + \Li\left(\frac 1z\right) = -\frac12 \log^2(-z) - \frac{\pi^2}6.$$
This identity is verified by differentiating both sides.  As $\Li' = \li'$ and $\Li(z) = \li(z)$ on the negative real axis, the identity holds for $\li$; hence, 
\begin{equation} \label{eq: identity}
\li(z) + \li\left(\frac1z\right) = -\frac12\log^2(-z) - \frac{\pi^2}6.
\end{equation}
(The branch of logarithm being used should be clear from context.)

\subsection{Deriving the formula for $F_{2,1}(x)$}

For a fixed positive value of $x$, we define the map $H_x: [0, \log\coth(x/2)] \to \bc$ as follows:
$$H_x(r) = \li(-e^{-r}\coth(x/2))-\li(e^{-r}\coth(x/2))+\li(-e^{-r}\tanh(x/2))-\li(e^{-r}\tanh(x/2)).$$

\begin{Lem}
$$\frac{d(\Re H_x)}{dr} = \log\frac{\coth x + \cosh r}{\coth x - \cosh r}$$
\end{Lem}

\begin{proof}
We first calculate $H_x'$ and then take real parts.  Given the definition of the dilogarithm and the fact that $\coth(x/2)+\tanh(x/2) = 2\coth x$, we have that
\begin{align*}
H_x'(r) &= \log(1+e^{-r}\coth(x/2))-\log(1-e^{-r}\coth(x/2))+\\
	&+\log(1+e^{-r}\tanh(x/2))-\log(1-e^{-r}\tanh(x/2))\\
	& = \resizebox{.9\hsize}{!}{$\log[(1+e^{-r}\coth(x/2))(1+e^{-r}\tanh(x/2))]-\log[(1-e^{-r}\coth(x/2))(1-e^{-r}\tanh(x/2))]$}\\
	& = \log[2e^{-r}(\cosh r + \coth x)] - \log[2e^{-r}(\cosh r - \coth x)]\\
	&= \log(\coth x + \cosh r) - \log(\cosh r - \coth x)\\
	& = \log \frac{\coth x + \cosh r}{\coth x - \cosh r} - i\pi.
\end{align*}
Given the domain for $H_x$, the argument of the logarithm above is always a positive real number.
\end{proof}

We therefore see that $F_{2,1}(x) = 2\cdot\Re[H_x(\log\coth(x/2)) - H_x(0)]$ as $\Omega_0= 2$. For a surface $S$ with boundary, let $\ell(\partial S)$ be the sum of the lengths of the boundary components.  Given the above we can now prove the following:

\medskip
\noindent
{\bf Corollary \ref{surface}.} {\it
Let $S$ be a compact hyperbolic surface with nonempty totally geodesic boundary.  Then
$$A_1(S) = \frac2{\ell(\partial S)} \sum_{\ell\in |O_S|}\left[\Li\left(-\tanh^2\frac {\ell}2\right)-\Li\left(\tanh^2\frac {\ell}2\right)+\frac {\pi^2}4\right].$$}

\begin{proof}
From the above formulation of $F_{2,1}(x)$, we get the following:
\begin{align*}
F_{2,1}(x) = 2\cdot\Re&\left[\li(a) + \li\left(\frac1a\right) - \li\left(-a\right)-\li\left(-\frac1a\right)+ \li\left(-\frac1{a^2}\right)-\li\left(\frac1{a^2}\right) - \frac {\pi^2}{4}\right],
\end{align*}
where $a = \coth\left(\frac x2\right)$.  From applying $\eqref{eq: identity}$ twice we see that:
$$\li(a) + \li\left(\frac1a\right) - \li\left(-a\right)-\li\left(-\frac1a\right) = \frac12(\log^2(a) - \log^2(-a)).$$
Recalling that $\log(-1) = i\pi$, for $a>0$ we have $\log(a) - \log(-a) = -\log(-1)$, so that simplifying the above
$$\li(a) + \li\left(\frac1a\right) - \li\left(-a\right)-\li\left(-\frac1a\right) = \frac{\pi^2}2 - i\pi\log a  \text{  for } a>0.$$
In particular, for positive values of $a$, the real part is always $\pi^2/2$.  As $\ell_i$ is always positive this identity holds in the case of interest. Furthermore, $\Li\left(\pm\tanh(\ell_i/2)\right) = \li\left(\pm\tanh(\ell_i/2)\right)$ as $\pm\tanh(\ell_i/2) \in (-1,1)$; hence, the result follows.
\end{proof}

\subsection{Asymptotics for $F_{2,1}(x)$}

We will use the following notation throughout the rest of the paper:  For $f,g:\br\to \br$ we will write $f\sim g$ if $\lim_{x\to\infty} [f(x)/g(x)] = 1$. This is clearly an equivalence relation on real-valued functions. Below we find the asymptotic behavior of $F_{2,1}(x)$ from our above result; we note that we  will also come to the same result later in the paper when we study the asymptotics of $F_{n,k}$ from the integral definition.

\begin{Prop}
Let $F_{2,1}(x)$ be defined as above, then $F_{2,1}(x)\sim 8xe^{-x}.$
\end{Prop}

\begin{proof}
We start with the following observation:
$$\lim_{x\to1} \frac{\Li(-x)-\Li(x)+\pi^2/4}{(x-1)\log(1-x)} = 1,$$
which is a direct application of L'H\^opital's rule and the definition of the dilogarithm.
From this, we can gather the following:
\begin{align*}
F_{2,1}(x) &\sim 2(\tanh^2(x/2)-1)\log(1-\tanh^2(x/2)) =4 \sech^2(x/2)\log\cosh(x/2)\\
&=4\left(\frac{2}{e^\frac x2+e^{-\frac x2}}\right)^2\log\left[e^\frac x2\left(\frac{1+e^{-x}}2\right)\right] \sim 8xe^{-x}
\end{align*}
\end{proof}


\section{Connecting with Bridgeman's identity in dimension 2}
\subsection{Liouville measure and Bridgeman's identity}
We first need to recall Bridgeman's identity.  Denote the space of oriented geodesics in $\bh^2$ by $G(\bh^2)$, then by identifying a geodesic with its endpoints in $\bS_\infty^{1}$ we see 
$$G(\bh^2) \cong (\bS_\infty^{1}\times \bS_\infty^{1})-\Delta,$$
where $\Delta$ denotes the diagonal and we view the geodesic $[x,y]\in G(\bh^2)$ as oriented from $x$ to $y$.  The \emph{ Liouville measure } $\mu$ is an $\mathrm{Isom}(\bh^2)$ invariant measure on $G(\bh^2)$.  If we identify $\bS_\infty^1$ with $\overline\br$, then $\mu$ is characterized by 
$$\mu((a,b)\times (c,d) ) = 2\left| \log\left| \frac{(a-c)(b-d)}{(a-d)(b-c)}\right|\right|$$
for disjoint intervals $(a,b)$ and $(c,d)$ in $\overline{\br}$.
The Liouville measure on $G(\bh^2)$ is derived from the Liouville measure on the unit tangent bundle $T_1(\bh^2)$, which is given by $dVd\theta$, where $dV$ is the hyperbolic volume measure and $d\theta$ is the standard measure on $\bS^1$ (see \cite{bonahon,nicholls} for details).  In fact, the natural fibering $T_1(\bh^2)\to G(\bh^2)$ is such that the volume measure on $T_1(\bh^2)$ is $d\mu\, dl$, where $dl$ is the length along a fiber. Note that the factor of 2 appears above so that $d\mu\,dl$ agrees with $dV\,d\theta$.

There are local coordinates for $G(\bh^2)$ where the Liouville measure can be written as a product measure.  Let $\gamma \in G(\bh^2)$ and $p\in \gamma$ a base point.  Let $\eta$ be a  geodesic intersecting $\gamma$, then $\eta$ is determined by the signed hyperbolic distance $s = \pm d(\gamma\cap \eta, p)$ coming from the orientation of $\gamma$ and the angle $\theta$ between $\gamma$ and $\eta$ measured from $\eta$ to $\gamma$.  In these local coordinates, we have
\begin{equation}\label{eq: measure}
d\mu_{(s,\theta)} = \sin\theta\, ds\, d\theta.
\end{equation}
These coordinates are described in the  appendices of \cite{bonahon}.

Given a hyperbolic surface $S$ with totally geodesic boundary, identify the universal cover of $S$ with a convex subset $\tilde S \subseteq \bh^2$.  Define $G(\tilde S)\subseteq G(\bh^2)$ to be the set of all geodesics intersecting $\tilde S$.  Let $\Gamma\subset \mathrm{Isom}(\bh^2)$ such that $S = \tilde S / \Gamma$, then we can set $G(S) = G(\tilde S)/\Gamma$ to be the space of geodesics in $S$.  The Liouville measure descends to a measure on $G(S)$.

Let $S$ be an orientable compact hyperbolic surface with nonempty totally geodesic boundary and given $v\in T_1(S)$ let $\al_v\in G(S)$ such that $\al_v'(t) = v$ for some $t\in \br$. For each orthogeodesic $\al_i$ set 
$W_i = \{v \in T_1(S)\colon \al_v \simeq \al_i \text{ rel } \partial S\},$
where $\simeq$ denotes homotopy equivalence in $S$. We then have $\Vol(T_1(S)) = 4\pi^2 |\chi(S)| = \sum \Vol(W_i)$.  If we define $L_S: G(S) \to \br$ by $L_S(g) = length(g)$, where length is measured in $S$, and set $\overline W_i = \{\al_v \in G(S)\colon v\in W_i\}$, then it was proved in \cite{b} that 
$$\Vol(W_i) = \int_{\overline W_i} L_S \,d\mu = 8\mathscr{L}\left(\sech^2 \frac {\ell_i}2\right),$$
where $\mathscr{L}(x) = \Li(x) + \frac12 \log|x|\log(1-x)$ for $x\leq 1$ is the Rogers dilogarithm.  This gives Bridgeman's identity:
$$\sum_i \mathscr{L}\left(\sech^2\frac{\ell_i}2\right) = \frac{\pi^2}2 |\chi(S)|.$$

\subsection{Random Variables}

Let $S$ be an orientable compact hyperbolic surface with nonempty totally geodesic boundary.  Given an angle $\theta\in (0,\pi)$,   we define a new random variable on $(\partial S, dm)$ as follows:  let $\gamma$ be a component of $\partial S$  and $x\in \gamma$.  Let $v\in T_x(S)$ such that the unit speed geodesic $\al_v$ resulting from the geodesic flow in the direction of $v$ makes an angle $\theta$ with $\gamma$ when measured from $\al_v$ to $\gamma$ (where the orientation of $\gamma$ is such that the surface is on the right).  Define $\al_x^\theta = \al_v$, then the function $L_\theta: \partial S \to \br$ defined by $L_\theta(x) = length(\al^\theta_x)$ is a random variable on $(\partial S, dm)$. We define its moments
$$A_k^\theta(S) = E[L_\theta^k] = \int_{\partial S} L_\theta^k\, dm.$$

As above, we can decompose our boundary into intervals associated to orthogeodesics: for $\al_i \in O_S$ let 
$$U^\theta_i = \{x\in \partial S\colon \al_x^\theta \simeq \al_i \text{ rel } \partial S\}.$$

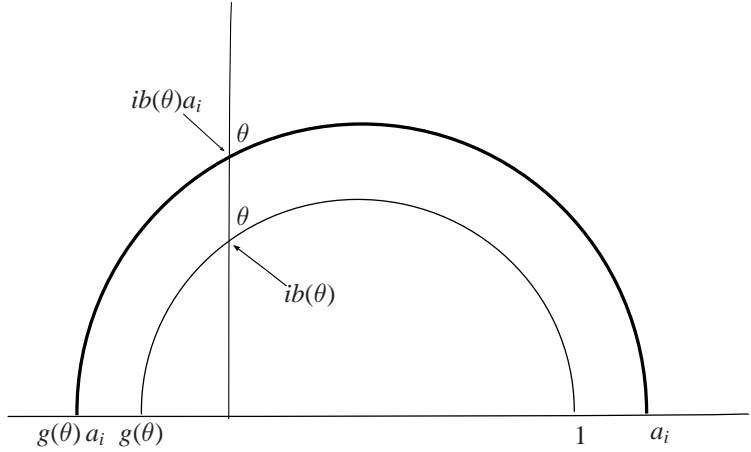
\begin{figure}[ht!]
\labellist 
\small\hair 2pt
\pinlabel $\theta$ at 148 180
\pinlabel $\theta$ at 148 128
\pinlabel $g(\theta)\,a_i$ at 40 -10
\pinlabel $g(\theta)$ at 84 -10
\pinlabel $1$ at 360 -10
\pinlabel $a_i$ at 410 -10
\pinlabel $ib(\theta)$ at 190 78
\pinlabel {$ib(\theta)a_i$} at 100 200
\endlabellist
\centering
\includegraphics[scale=0.6, clip=false]{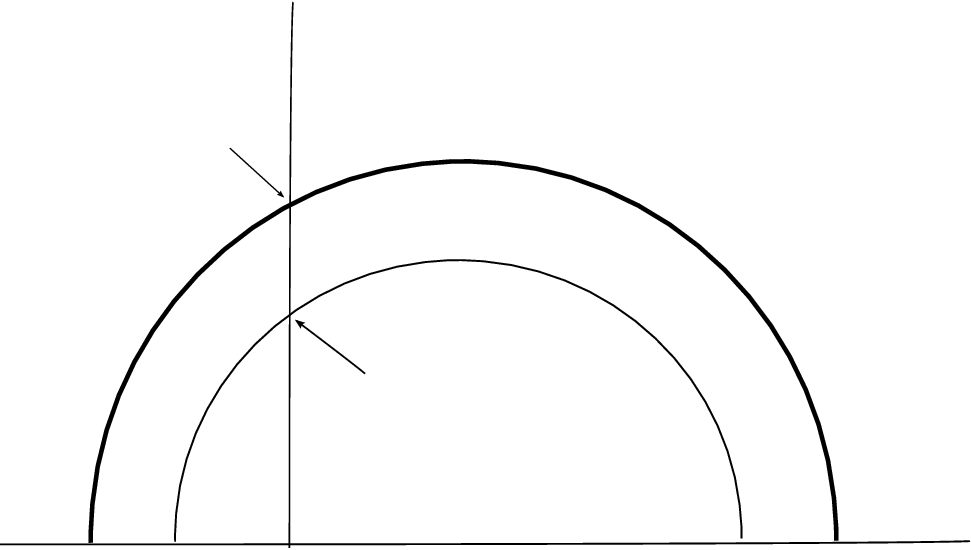}
\caption{The figure shows the $\theta$-projection of the geodesic $[1,a_i]$ to $[0,\infty]$ as in Lemma \ref{lem: theta}.}
\label{fig:figure2}
\end{figure}

\begin{Lem} \label{lem: theta}
For every $\theta \in (0,\pi)$, the set $U_i^\theta$ is an interval of length $2\log\coth(\ell_i/2)$.
\end{Lem}

\begin{proof}
Let $\gamma_1, \gamma_2$ be the components of $\partial S$ such that $\al_i$ travels from $\gamma_1$ to $\gamma_2$ (possibly $\gamma_1=\gamma_2$).  We may then put this picture in the upper half plane with $\tilde{\gamma}_1 = [0,\infty]$ and $\tilde\gamma_2=[1,a_i]$, where $a_i= \coth^2(\ell_i/2)$. For $x\in U_i^\theta$ there exists a unique lift $\tilde\al_x^\theta$ intersecting $\tilde \gamma_1$ in angle $\theta$ and also intersecting $\tilde \gamma_2$.  As in the proof of Basmajian's identity, we see $U_i^\theta$ lifts to the  $\theta$-projection of $\gamma_2$ onto $\gamma_1$.  Define $g(\theta)$ such that the geodesic $\be=[g(\theta),1]$ intersects $\tilde \gamma_1$ at angle $\theta$ as in Figure \ref{fig:figure2}. Define $b(\theta)$ so that $\beta$ intersects $\tilde \gamma_1$ at the point $ib(\theta)$.  
Observe that the geodesic $[a_ig(\theta), a_i]$ intersects $\tilde\gamma_1$ at angle $\theta$ at the point $ia_ib(\theta)$ as it is the image of $\be$ under the M\"obius transformation $z\mapsto a_iz$.  In particular, the length of $U_i^\theta$ is $\log(a_i)$.
\end{proof}

In the above proof we see that the picture of $U_i^\theta$ only depends on $\ell_i$, so as in Theorem \ref{ortho} we have that there exists functions 
$$F_k^\theta(\ell_i) = \int_{U_i^\theta} L_\theta^k \,dm,$$
such that
$$A_k^\theta(S) = \sum_i F_k^\theta(\ell_i).$$
In particular, each $A_k^\theta$ gives an orthospectrum identity.

\begin{Prop}
For  $\theta\in(0,\pi)$, the identity for $A_0^\theta(S)$ is Basmajian's identity.
\end{Prop}

\begin{proof}
As $m(\bigcup_i U_i^\theta) =1$ and $U_i^\theta\cap U_j^\theta = \emptyset$ for $i\neq j$, we have 
$$1 = \sum_i m(U_i^\theta) = \frac{1}{\ell(\partial S)}\sum_i 2\log\coth(\ell_i/2).$$
\end{proof}

We now have the following connection between $A_1(S)$ and Bridgeman's identity:

\begin{Thm} Suppose $S$ is a compact hyperbolic surface with nonempty totally geodesic boundary.
Let $F_S: [0,\pi] \to \br$  be defined by 
$$F_S(\vp)=\ell(\partial S) \int_0^\vp A_1^\theta (S) \sin\theta \, d\theta,$$
then the identity for $F_S(\pi)$ is Bridgeman's identity and $F'(\pi/2) = A_1(S)$.
\end{Thm}

\begin{proof}
Let $\tilde S$ be the universal cover of $S$ identified with a convex subset of the upper half plane $\bh$.  Let $\al_i$ be an orthogeodesic traveling from the boundary component $\gamma_1$ to the component $\gamma_2$.  Assume that the geodesic $[0, \infty]\subset \bh^2$ is a lift of $\gamma_1$ and the geodesic $[1, a_i]\subset \bh^2$ is a lift of $\gamma_2$, where $a_i = \log\coth(\ell_i/2)$.  As in the proof of Lemma \ref{lem: theta}, we lift the set $\overline W_i$ to the set $\wt W_i = (-\infty,0)\times (1, a_i)\subset G(\tilde S)$.  Then every geodesic $[x,y]\in \wt W_i$  can be given coordinates $(s,\theta)$, where
$$ s([x,y]) = \log \left(\frac{[0,\infty]\cap [x,y]}{i}\right)$$
 and $\theta([x,y])$ is the the angle from $[x,y]$ to $[0,\infty]$.  Using these local coordinates the Liouville measure can be written as in \eqref{eq: measure} and from the notation above it follows that 
$$\Vol(W_i) = \int_{\overline{W}_i} L_S\, d\mu = \int_0^\pi \int_{\log(b(\theta))}^{\log(b(\theta) a_i)} L_\theta(s)\sin\theta \, ds\,d\theta $$
and thus
\begin{align*}
F_S(\pi) &= \ell(\partial S)\int_0^\pi A_1^\theta(S)\sin\theta\, d\theta\\
& =\sum_i \int_0^\pi \int_{\log(b(\theta))}^{\log(b(\theta) a_i)} L_\theta(s)\sin\theta \, ds\,d\theta\\
& = \sum_i \Vol(W_i)  \\
&=8\sum_i\mathscr{L}\left(\sech^2\frac{\ell_i}2\right).
\end{align*}

Using the notation from the beginning of the section, we see that almost every element of $G(S)$ can be realized as $\al_p^\theta$ for $p\in \partial S$ and $\theta\in (0,\pi)$ implying that $\partial S\times (0,\pi)$ is full measure in $G(S)$.  In particular,
$$F_S(\pi) = \ell(\partial S) \int_0^\pi A_1^\theta(S) \sin\theta \, d\theta = \int_0^\pi \int_{\partial S} L_\theta \sin\theta\, ds\,d\theta = \int_{G(S)}L_S\, d\mu = 4\pi^2|\chi(S)|.$$
This shows that the identity for $F_S(\pi)$ is Bridgeman's identity.  Further, it is clear from the definition that $F_S'(\pi/2) = A_1(S)$.
\end{proof}
\noindent
{\bf Remark.}
This also shows that $A_1^\theta(S)$ is finite.  It can also be seen that $A_k^\theta(S)$ is finite for all $k$ using the same approach as in the proof for the finiteness of $A_n(M)$ given earlier.

\section{Odd Dimensions}

In this section we will write an explicit formula for $A_1(M^n)$ with $n$ odd.  For $n$ odd, we can simplify the integral in the definition of $F_{n,k}$ by using the substitution $u=\cosh r$ to get:
$$F_{n,1}(x) = \Omega_{n-2}\int_1^{\coth x} (u^2-1)^{\frac{n-3}2}\log\frac{\coth x + u}{\coth x - u} \, du.$$
An elementary calculation gives the following integrals (up to a constant) when $m$ is even:
$$\resizebox{\hsize}{!}
{$F^+_m(u,y) = \int u^m \log(y +  u) \, du = \frac{1}{m+1}\left[\left(u^{m+1}+y^{m+1}\right)\log(y+u) + \sum_{k=1}^{m+1}(-1)^{m-k}\frac{y^{m-k+1}u^k}{k}\right]$
}$$
$$F^-_m(u,y) = \int u^m\log(y- u)\, du = \frac{1}{m+1}\left[\left(u^{m+1}-y^{m+1}\right)\log(y-u) - \sum_{k=1}^{m+1}\frac{y^{m-k+1}u^k}{k}\right].$$
\noindent
Now set 
$$f_m(x) = F^+_m(\coth x, \coth x)-F^+_m(1,\coth x)+F^-_m(1,\coth x) - \lim_{u\to (\coth x)^-} F^-_m(u,\coth x).$$
After some routine manipulation and simplification, we find:
$$f_m(x) = \frac{2\coth^{m+1}(x)}{m+1}\left[\log(2\cosh x)-x\tanh^{m+1}(x)+\sum_{k=1}^{\frac m2} \frac{1-\tanh^{2k}(x)}{2k}\right].$$
If we expand out the binomial in $F_{n,1}(x)$, we find
$$F_{n,1}(x) = \Omega_{n-2}\sum_{j=0}^{\frac{n-3}2}(-1)^{\frac{n-3}2-j}\binom{\frac{n-3}2}{j}f_{2j}(x).$$
We then immediately have:

\medskip
\noindent
{\bf Corollary \ref{odd}.} {\it
Let $M$ be an $n$-dimensional compact hyperbolic manifold with nonempty totally geodesic boundary where $n$ is odd.  Then
\begin{align*}
\resizebox{\hsize}{!}{$\displaystyle{A_1(M) = \frac{2\Omega_{n-2}}{Vol(\partial M)}\sum_{\ell_i\in |O_M|}\sum_{j=0}^{\frac{n-3}2}\frac{(-1)^{\frac{n-3}2-j}\binom{\frac{n-3}2}{j}}{2j+1}\coth^{2j+1}(\ell_i)\left[\log(2\cosh\ell_i)-\ell_i\tanh^{2j+1}(\ell_i)+\sum_{k=1}^j\frac{1-\tanh^{2k}(\ell_i)}{2k}\right]}.$}
\end{align*}}


\section{The Asymptotics of $F_{n,k}$}

In this section, we explore the asymptotic behavior of the $F_{n,k}$'s. 

\begin{Thm}
For all $n,k\in \bz^+$ $$\lim_{x\to\infty} \frac{e^{-(n-1)x}F_{n,k}(x)}{x^k} = \frac{2^{n+k-1}\Omega_{n-2}}{n-1},$$ 
Furthermore, for $n$ odd we have 
$$\lim_{x\to0}{x^{n-2}}F_{n,1}(x) = \frac{2}{n-2}[\log 2 + \frac12 H_{(n-1)/2}],$$
where $H_n$ is the $n^{th}$ harmonic number.
\end{Thm}

\begin{proof}
Recall that $F_{n,k}(x) = \Omega_{n-2} \int_0^{\log\coth(x/2)} L^k_x(r) \sinh^{n-2}(r)\, dr$.  Using the substitution $u=\cosh r$, we have 
$$F_{n,k}(x) = \Omega_{n-2} \int_1^{\coth x} (u^2-1)^{\frac{n-3}2} \left[\log\left(\frac{\coth x+u}{\coth x- u}\right)\right]^k du.$$
For the moment, let $n\geq 3$, so that $(n-3)/2 \geq 0$, then
\begin{align*}
F_{n,k}(x)& \geq 2^{\frac{n-3}2}\Omega_{n-2}\int_1^{\coth x} (u-1)^{\frac{n-3}2}\left[\log\left(\frac{\coth x + u}{\coth x-u}\right)\right]^k du\\
F_{n,k}(x) & \leq (\coth x +1)^{\frac{n-3}2}\Omega_{n-2}\int_1^{\coth x} (u-1)^{\frac{n-3}2}\left[\log\left(\frac{\coth x + u}{\coth x-u}\right)\right]^k du.
\end{align*}
As $\coth x+1$ approaches 2 in the limit as $x$ goes to infinity, we see from the above two inequalities that 
$$F_{n,k}(x) \sim 2^{\frac{n-3}2}\Omega_{n-1}\int_1^{\coth x} (u-1)^{\frac{n-3}2}\left[\log\left(\frac{\coth x + u}{\coth x-u}\right)\right]^k du.$$
In the case $n=2$, the inequalities above are reversed, but yield the same result; hence, what follows will hold for all $n$.  We now look at the following two inequalities:
$$\resizebox{\hsize}{!}{$
\int_1^{\coth x} (u-1)^{\frac{n-3}2}\left[\log\left(\frac{\coth x + u}{\coth x-u}\right)\right]^k du  \geq \int_1^{\coth x} (u-1)^{\frac{n-3}2}\left[\log(\coth x+1) - \log(\coth x -u)\right]^k du$}$$
$$\resizebox{\hsize}{!}{$
\int_1^{\coth x} (u-1)^{\frac{n-3}2}\left[\log\left(\frac{\coth x + u}{\coth x-u}\right)\right]^k du  \leq \int_1^{\coth x} (u-1)^{\frac{n-3}2}\left[\log(2\coth x) - \log(\coth x -u)\right]^k du$}$$

Note for large $x$ that $\coth x - u<1$ for all $u\in [1,\coth x]$, so that $\log(\coth x - u) < 0$. As both $\log(2\coth x)$ and $\log(\coth x+1)$ limit to $\log 2$, we see that both the integrals  in the inequalities are asymptotic to $\int_1^{\coth x} (u-1)^{(n-3)/2} [\log(\coth x- u)]^k du$.  Let us write $a(x) = \coth x -1$ and $v = \frac{u-1}a$, so that we now have
$$\resizebox{\hsize}{!}{$
F_{n,k}(x) \sim  2^{\frac{n-3}2}\Omega_{n-2}a^{\frac {n-1}2}\int_0^1 v^{\frac{n-3}2} [\log(a-av)]^k dv =  2^{\frac{n-3}2}\Omega_{n-2}a^{\frac {n-1}2}\int_0^1 v^{\frac{n-2}2} \left[\log a + \log\left(1-v\right)\right]^k dv.$}$$
As $\int_0^1 v^{(n-3)/2} [\log(1-v)]^m dv$ is finite for all $m$, we find that 
$$F_{n,k}(x) \sim (-1)^k 2^{\frac{n-3}2} \Omega_{n-2}(\log a)^{k}a^{\frac {n-1}2} \int_0^1 v^{\frac{n-3}2}  dv =\frac {2^\frac {n-1}2 \Omega_{n-2}}{n-1} \left(\log\frac1a\right)^k a^{\frac{n-1}2}.$$
Since, $a(x) = \coth x -1 \sim 2e^{-2x}$, we get the stated result.

When $n$ is odd, we have the following when $x$ approaches 0: As $x$ tends to 0, it is easy to see that $\tanh^{m+1}(x)f_m(x)$ is finite.  As $\lim_{x\to 0}[x \coth x]$ is finite, we see that $\lim_{x\to 0} {x^{m+1}}f_m(x)<\infty$.  Again, as $F_{n,1}(x)$ is a sum of the $f_m$'s, the largest exponent dominates, which gives the result.
\end{proof}

\section{The Moment Generating Function in Dimension 3}

Let $M$ be a hyperbolic 3-manifold with totally geodesic boundary and let $S = \partial M$.  We define the {\it moment-generating function} $M_L(t) = E[e^{tL}]$, where $E[X]$ denotes the expected value of a random variable $X$ with respect to our probability measure $dm = dV/V(\partial M)$.  The moment-generating function encodes all the moments of $L$ in its derivatives:  $A_k(M)=E[L^k] = M_L^{(k)}(0)$.  In particular, by calculating $M_L(0)$ we will recover Basmajian's identity and $A_1(M)$ by calculating $M_L'(0)$. The goal of this section is to prove that following theorem:

\begin{Thm}\label{generating}
Let $M$ be a compact hyperbolic 3-manifold with totally geodesic boundary $S$ and let $\delta$ be the Hausdorff dimension of the limit set of $M$.  For $t \in (\delta-2, 2-\delta)$
$$M_L(t) = \frac{4\pi}{V(S)}\sum_{\ell_i\in |O_M|} \coth(\ell_i)\cdot B\left(\frac{1-\tanh \ell_i}{2},1-t,1+t\right),$$
where $B$ is the incomplete beta function.
\end{Thm}

\subsection{Hypergeometric Function and Incomplete Beta Function}

The {\it hypergeometric functions} $_2F_1(a,b,c,z)$ for $z\in\bc$ with $|z|<1$ are given by the power series:

$$_2F_1(a,b,c,z) = \sum_{n=0}^\infty \frac{(a)_n(b)_n}{(c)_n} \frac{z^n}{n!},$$
provided $c\notin \bz^{\leq 0}$ and where
$$(a)_n = \left\{ 
\begin{array}{lr}
1 & \text{for }n = 0 \\
a(a+1)\cdots(a+n-1) & \text{for }n >0
\end{array} \right. .$$
We will use the following identity below:
$$(1-z)^{-a} = {_2F_1}(a,1,1,z).$$

The {\it incomplete beta functions} $B(x,a,b)$ are defined as
$$B(x,a,b) = \int_0^x s^{a-1}(1-s)^{b-1}ds.$$
We can also write an incomplete beta function in terms of a hypergeometric function as follows (see \cite{luke}):
\begin{equation} \label{eqn: beta}
B(x,a,b) = \frac{x^a}a {_2F_1}(a,1-b,a+1,x).
\end{equation}

We present two calculations as technical lemmas that will allow us to derive the moment generating function.

\begin{Lem}
$$\frac{\partial}{\partial x} {_2F_1}(1+t,t,2+t,x) = \frac{1+t}x \left[(1-x)^{-t} - {_2F_1}(1+t,t,2+t,x)\right]$$
\end{Lem}

\begin{proof}
We calculate:
\begin{align*}
\frac{\partial}{\partial x} {_2F_1}(1+t,t,2+t,x) &= \sum_{n=1}^\infty \frac{(t)_n(1+t)_n}{(2+t)_n} \frac{x^{n-1}}{(n-1)!} \\
&= \frac{1+t}x \sum_{n=1}^\infty(t)_{n}\frac{n}{t+n+1}\frac{x^n}{n!}\\
&= \frac{1+t}x \sum_{n=1}^\infty \left[(t)_n-\frac{(t)_n(1+t)}{t+n+1}\right]\frac{x^n}{n!} \\
&= \frac{1+t}x\left[ \sum_{n=0}^\infty (t)_n\frac{x^n}{n!}- \sum_{n=0}^\infty \frac{(t)_n(1+t)_n}{(t+2)_n} \frac{x^n}{n!}\right] \\
& = \frac{1+t}{x}[{_2F_1}(t,1,1,x)-{_2F_1}(1+t,t,2+t,x)]\\
&=\frac{1+t}x \left[(1-x)^{-t} - {_2F_1}(1+t,t,2+t,x)\right]
\end{align*}
\end{proof}

\begin{Lem}
Let $g(u,a,t) = (1+t)^{-1}(a+u)^{t+1}(2a)^{-t}{_2F_1}\left(1+t,t,2+t, \frac{a+u}{2a}\right),$ then
$$\frac{\partial g}{\partial u} = \left(\frac{a+u}{a-u}\right)^t.$$
\end{Lem}

\begin{proof}
This is an immediate consequence of the previous lemma.
\end{proof}

\subsection{Proof of Theorem \ref{generating}}
We can now find the moment generating function of $L$.  
\begin{proof}[Proof of Theorem \ref{generating}]
Let $S = \partial M$ and recall that $\Omega_1 = 2\pi$.  By definition, 
\begin{align*}
M_L(t) &= E[e^{tL}] = \int_S e^{tL}dm = \sum_i \int_{U_i} e^{tL_i}dm \\ 
&= \frac{2\pi}{V(S)}\sum_i\int_{0}^{\log(\coth(\ell_i/2))} \left(\frac{\coth(\ell_i)+\cosh r}{\coth(\ell_i)-\cosh r}\right)^t \sinh r\, dr.\\
&= \frac{2\pi}{V(S)}\sum_i\int_1^{\coth(\ell_i)}\left(\frac{\coth(\ell_i)+u}{\coth(\ell_i)-u}\right)^t du,
\end{align*}
where $u = \cosh r$. 
By replacing $L^k$ with $e^{tL}$ in the proof of Theorem \ref{finiteness}, we see that $\int_S e^{tL}dm$ converges for $t < 2-\delta$. Note that if $t\in (2-\delta, \delta-2)$ then $|t| < 1$.  
From the above lemma, we then have that 
$$M_L(t) = (2\pi/V(S))\sum_{\ell_i\in |O_M|}[g(\coth(\ell_i),\coth(\ell_i),t) - g(0,\coth(\ell_i),t)].$$ 
After expanding the above terms using the definition of $g$, some simplifications get us to:
$$\resizebox{\hsize}{!}{$\displaystyle{
M_L(t) =\frac{2\pi}{V(S)} \sum_{\ell_i\in |O_M|} \frac{\coth(\ell_i)}{1+t}\left[2{_2F_1}(1+t,t,2+t,1)-\frac12\left(\frac{1+\tanh(\ell_i)}{2}\right)^{t+1}{_2F_1}\left(1+t,t,2+t,\frac{1+\tanh(\ell_i)}2\right)\right].}$}$$
By \eqref{eqn: beta} this becomes
$$M_L(t) =\frac{4\pi}{V(S)} \sum_{\ell_i\in |O_M|} \coth(\ell_i)\left[B(1,1+t,1-t)-B\left(\frac{1+\tanh(\ell_i)}2,1+t,1-t\right)\right].$$
It is left to investigate $B(1,1+t,1-t) - B(a, 1+t,1-t)$:
$$\resizebox{\hsize}{!}{$
B(1,1+t,1-t)-B(a,1+t,1-t) = \int_a^1s^t(1-s)^{-t}ds = -\int_{1-a}^0(1-u)^tu^{-t} du = B(1-a,1-t,1+t),$}$$
where $u = 1-a$.  Therefore, we can conclude
$$M_L(t) = \frac{4\pi}{V(S)}\sum_{\ell_i\in |O_M|}\coth(\ell_i) \cdot B\left(\frac{1-\tanh(\ell_i)}2,1-t,1+t\right).$$
\end{proof}

\subsection{Recovering Basmajian's Identity in Dimension 3}

As $M_L(0) = 1$ we have 
$$1 = \frac{4\pi}{V(S)} \sum_{\ell_i\in |O_M|} \coth(\ell_i)\cdot B\left(\frac{1-\tanh(\ell_i)}2,1,1\right)$$
and as $B(a,1,1) = a$, we have
$$V(S) = \sum_{\ell_i\in |O_M|} 2\pi(\coth(\ell_i) - 1) = \sum_{\ell_i\in |O_M|}\frac{2\pi e^{-\ell_i}}{\sinh(\ell_i)} = \sum_{\ell_i\in |O_M|} V_2(\log(\coth(\ell_i/2))),$$
where $V_2(r)$ is the area of a hyperbolic circle of radius $r$.

\footnotesize
\nocite{*}
\bibliographystyle{gtart}	
\bibliography{average-length}

\end{document}